\documentclass[12pt,reqno]{amsart}

\usepackage{amsfonts, amsthm, amsmath}
\allowdisplaybreaks[4]

\usepackage{rotating}

\usepackage{tikz}

\usepackage{graphics}

\usepackage{amssymb}

\usepackage{amscd}

\usepackage{graphicx}

\usepackage{epstopdf}

\usepackage[latin2]{inputenc}

\usepackage{t1enc}

\usepackage[mathscr]{eucal}

\usepackage{indentfirst}

\usepackage{graphicx}

\usepackage{graphics}

\usepackage{pict2e}

\usepackage{mathrsfs}

\usepackage{enumerate}
\usepackage[pagebackref]{hyperref}
\hypersetup{backref, pagebackref, colorlinks=true}
\usepackage{cite}
\usepackage{color}
\usepackage{epic}
\usepackage{hyperref} 
\numberwithin{equation}{section}
\theoremstyle{plain}

\newtheorem{theorem}{Theorem}[section]

\theoremstyle{definition}

\newtheorem{remark}[theorem]{Remark}

\newtheorem{?}[theorem]{Problem}

\def\boxit#1{\leavevmode\hbox{\vrule\vtop{\vbox{\kern.33333pt\hrule
    \kern1pt\hbox{\kern1pt\vbox{#1}\kern1pt}}\kern1pt\hrule}\vrule}}

\newcommand{\bl}[1]{\uppercase\expandafter{\romannumeral#1}}

\usepackage{collectbox}



\begin{document}

\title[Multi-dimensional $q$-summations]{Multi-dimensional $q$-summations and multi-colored partitions}

\author[S. Chern]{Shane Chern}
\address[Shane Chern]{Department of Mathematics, The Pennsylvania State University, University Park, PA 16802, USA}
\email{shanechern@psu.edu}

\author[S. Fu]{Shishuo Fu}
\address[Shishuo Fu]{College of Mathematics and Statistics, Chongqing University, Huxi campus LD506, Chongqing 401331, P.R. China}
\email{fsshuo@cqu.edu.cn}

\author[D. Tang]{Dazhao Tang}

\address[Dazhao Tang]{College of Mathematics and Statistics, Chongqing University, Huxi campus LD208, Chongqing 401331, P.R. China}
\email{dazhaotang@sina.com}

\date{\today}

\begin{abstract}
Motivated by Alladi's recent multi-dimensional generalization of Sylvester's classical identity, we provide a simple combinatorial proof of an overpartition analogue, which contains extra parameters tracking the numbers of overlined parts of different colors. This new identity encompasses a handful of classical results as special cases, such as Cauchy's identity, and the product expressions of three classical theta functions studied by Gauss, Jacobi and Ramanujan.
\end{abstract}



\maketitle


\section{Introduction}\label{sec1}
In 1882, Sylvester \cite{Syl} discovered the following identity:
\begin{align}
\label{eq1.1}
(-aq;q)_{\infty}=1+\sum_{k=1}^{\infty}\frac{a^{k}q^{(3k^{2}-k)/2}(-aq;q)_{k-1}(1+aq^{2k})}{(q;q)_{k}}.
\end{align}
Here and in the sequel, we use the standard $q$-series notation \cite{Andr}:
\begin{align*}
(a;q)_n:=&\prod_{k= 0}^{n-1} (1-aq^k),\\
(a;q)_\infty:=&\prod_{k= 0}^\infty (1-aq^k).
\end{align*}
The case $a=-1$ in \eqref{eq1.1} yields Euler's celebrated pentagonal number theorem:
\begin{align}\label{theta:3}
(q;q)_{\infty} &=1+\sum_{k=1}^{\infty}(-1)^{k}q^{(3k^{2}-k)/2}(1+q^{k})\notag\\
 &=\sum_{k=-\infty}^{\infty}(-1)^{k}q^{(3k^{2}-k)/2}.
\end{align}
The right-hand side of \eqref{theta:3} is one of the three classical theta functions studied by Gauss, Jacobi and Ramanujan. The other two allow similar product representations as follows.
\begin{align}
\label{theta:2}
\frac{(q;q)_{\infty}}{(-q;q)_{\infty}} &=\sum_{k=-\infty}^{\infty}(-1)^{k}q^{k^2},\\
\label{theta:4}
\frac{(q^2;q^2)_{\infty}}{(-q;q^2)_{\infty}} &=\sum_{k=-\infty}^{\infty}(-1)^{k}q^{2k^{2}-k}.
\end{align}
Empirically, properties enjoyed by one of these theta functions are usually shared by the other two, as witnessed by a recent work of the second and third authors \cite{FT2}. Our current investigation is of no exception (see Remark~\ref{rmk}).

A \emph{partition} of a nonnegative integer $n$ is a weakly decreasing sequence of positive integers whose sum equals $n$. Based on the observation that the left-hand side of \eqref{eq1.1} is the generating function of \emph{strict partitions} (i.e.~partitions into distinct parts), Sylvester proved his identity combinatorially by analyzing the Ferrers graphs of strict partitions in terms of their Durfee squares. The interested readers may refer to \cite{Andr3} for details.

In a recent paper \cite{All}, Alladi further considered $r$-colored strict partitions (i.e.~$r$ copies of strict partitions attached with colors $a_{1}$, $a_{2}$, $\ldots$, $a_{r}$). He then naturally generalized Sylvester's identity to a multi-dimensional summation, which can be stated as follows.

\begin{theorem}[Alladi]\label{thm:alladi}
We have
\begin{align}
&(-a_{1}q;q)_{\infty}(-a_{2}q;q)_{\infty}\cdots(-a_{r}q;q)_{\infty}\nonumber\\
&\quad=1+\sum_{N=1}^{\infty}q^{N^{2}}\prod_{j=1}^{r}(-a_{j}q;q)_{N-1} \sum_{i_{1}+i_{2}+\cdots+i_{r}=N}\frac{a_{1}^{i_{1}}a_{2}^{i_{2}}\cdots a_{r}^{i_{r}}q^{\binom{i_1}{2}+\binom{i_2}{2}+\cdots+\binom{i_r}{2}}}{(q;q)_{i_{1}}(q;q)_{i_{2}}\cdots(q;q)_{i_{r}}}\nonumber\\
&\quad\quad\quad\quad\quad\times \left(1+\sum_{s=1}^{r}q^{i_1+i_2+\cdots+i_s}a_s q^N \prod_{k=1}^{s-1}\left(1+a_k q^N\right)\right).\label{id:alladi}
\end{align}
\end{theorem}

We remark that Alladi's original identity (cf.~\cite[Eq.~(4.8)]{All}) involves some combinatorial statistics. However, he then showed in his equation (4.9) that the combinatorial statistics can be replaced and hence the multiple summation can be stated as above. In fact, he provided both an analytic and a combinatorial proof of \eqref{id:alladi}. Nonetheless, his combinatorial proof is complicated to some extent. This motivated us to give a simplified combinatorial proof. During the course, we are naturally led to the following {\em$r$-colored overpartition} (see Section 3 for the definition) analogue:
\begin{theorem}\label{thm:CFT}
We have
\begin{align}
&\frac{(-a_{1} z_1 q;q)_{\infty}(-a_{2} z_2 q;q)_{\infty}\cdots(-a_{r} z_r q;q)_{\infty}}{(a_{1}q;q)_{\infty}(a_{2}q;q)_{\infty}\cdots(a_{r}q;q)_{\infty}}\nonumber\\
&\quad=1+\sum_{N=1}^{\infty}q^{N^{2}}\prod_{j=1}^{r}\frac{(-a_{j} z_j q;q)_{N-1}}{(a_{j}q;q)_{N-1}} \sum_{i_{1}+i_{2}+\cdots+i_{r}=N}\frac{a_{1}^{i_{1}}a_{2}^{i_{2}}\cdots a_{r}^{i_{r}}(-z_1;q)_{i_{1}}(-z_2;q)_{i_{2}}\cdots(-z_r;q)_{i_{r}}}{(q;q)_{i_{1}}(q;q)_{i_{2}}\cdots(q;q)_{i_{r}}}\nonumber\\
&\quad\quad\quad\quad\quad\times \left(1+\sum_{s=1}^{r}q^{i_1+i_2+\cdots+i_{s-1}}\left(1+z_s q^{i_s}\right) \frac{a_s q^N}{1-a_s q^N}  \prod_{k=1}^{s-1} \frac{1+a_k z_k q^N}{1-a_k q^N}\right).\label{eq:over-gen}
\end{align}
\end{theorem}

The rest of this paper is organized as follows. In Section~\ref{sec:simple proof}, we provide a simplified combinatorial proof of \eqref{id:alladi}. In Section~\ref{sec:generalization}, we apply our approach to multi-colored overpartitions and prove \eqref{eq:over-gen}. We close with some remarks to motivate further investigations.

\section{A simple combinatorial proof of Theorem \ref{thm:alladi}}\label{sec:simple proof}
We could have proven Theorem~\ref{thm:CFT} directly and shown how to make appropriate substitutions for the variables to imply Theorem~\ref{thm:alladi}. However, we decide to warm the reader up by beginning with the proof of Theorem~\ref{thm:alladi}, since the combinatorial analysis in this case is simpler.

We first assume the following generalized order of parts in an $r$-colored (strict) partition:
$$1_{a_1}<1_{a_2}<\cdots<1_{a_r}<2_{a_1}<2_{a_2}<\cdots<2_{a_r}<3_{a_1}<\cdots.$$
When we plot the Ferrers graphs of these $r$-colored partitions, we color only the last node on the right of each row; the remaining nodes are uncolored.

\begin{figure}[ht]
\caption{Four blocks in a partition}\label{fig:block}
\vspace{1em}
\includegraphics{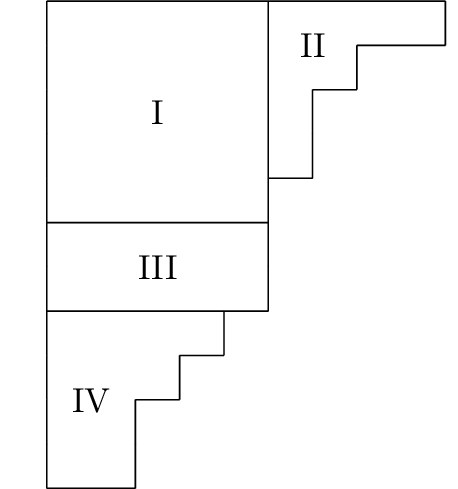}
\end{figure}

For an $r$-colored partition $\lambda$, its \emph{Durfee square} $D$ is defined to be the largest square of nodes contained within the Ferrers graph. We denote it as Block \uppercase\expandafter{\romannumeral1} in Fig.~\ref{fig:block}. We then denote by Block \uppercase\expandafter{\romannumeral2} the portion to the right of the Durfee square. Furthermore, the parts below the Durfee square that have the same size as the size of the Durfee square form Block \uppercase\expandafter{\romannumeral3}. At last, the portion below Block \uppercase\expandafter{\romannumeral3} is called Block \uppercase\expandafter{\romannumeral4}.

We remark that in Block \uppercase\expandafter{\romannumeral2} we also allow $0$ as a part. In this sense, we do not color any nodes in Block \uppercase\expandafter{\romannumeral1}, while instead we color the $0$ parts in Block \uppercase\expandafter{\romannumeral2}.

Now we are ready to write the generating function of each block combinatorially.

Let $N\ge 1$ be the size of the Durfee square $D$.

\emph{Block \uppercase\expandafter{\romannumeral1}}: Note that all nodes in $D$ are uncolored. Hence the generating function of $D$ is simply
\begin{align}
\label{cd1}
q^{N^{2}}.
\end{align}

\emph{Block \uppercase\expandafter{\romannumeral4}}: Note that Block \uppercase\expandafter{\romannumeral4} can be regarded as an $r$-colored strict partition with largest part $\le N-1$. Hence its generating function is
\begin{align}
\label{cd2}
\prod_{j=1}^{r}(-a_{j}q;q)_{N-1}.
\end{align}

\emph{Blocks \uppercase\expandafter{\romannumeral2} \& \uppercase\expandafter{\romannumeral3}}: We discuss the following two cases:
\begin{enumerate}[\indent 1).]
\item If Block \bl{3} is empty, then the generating function of Block \bl{2} is
\begin{align*}
\sum_{i_{1}+i_{2}+\cdots+i_{r}=N}\frac{a_{1}^{i_{1}}a_{2}^{i_{2}}\cdots a_{r}^{i_{r}}q^{\binom{i_1}{2}+\binom{i_2}{2}+\cdots+\binom{i_r}{2}}}{(q;q)_{i_{1}}(q;q)_{i_{2}}\cdots(q;q)_{i_{r}}}.
\end{align*}
\item If Block \bl{3} is not empty, then we assume that the part on the top of Block \bl{3} is colored by $a_s$ with $1\le s\le r$. Then the generating function of Block \bl{3} is given by
\begin{align*}
a_s q^N \prod_{k=1}^{s-1} \left(1+a_k q^N\right).
\end{align*}
Furthermore, in this case, we only allow $0$ colored by $a_{s+1}$, $\ldots$, $a_r$ as a part in Block \bl{2} to ensure that the whole is an $r$-colored strict partition. We assume that there are $i_t$ parts colored by $a_t$ in Block \bl{2} for each $1\le t\le r$. Then $i_{1}+i_{2}+\cdots+i_{r}=N$. For $1\le t_1\le s$, all distinct parts colored by $a_{t_1}$ can be regarded as a strict partition with exactly $i_{t_1}$ parts in the conventional sense (i.e.~$0$ is not allowed as a part), and hence have generating function
\begin{align*}
\frac{a_{t_1}^{i_{t_1}}q^{\binom{i_{t_1}}{2}+i_{t_1}}}{(q;q)_{i_{t_1}}}.
\end{align*}
For $s+1\le t_2\le r$, these parts colored by $a_{t_2}$ form a strict partition with either $i_{t_2}$ or $i_{t_2}-1$ parts in the conventional sense, which has generating function
\begin{align*}
\frac{a_{t_2}^{i_{t_2}}q^{\binom{i_{t_2}}{2}}}{(q;q)_{i_{t_2}}}.
\end{align*}
\end{enumerate}
We conclude that the generating function of Blocks \bl{2} \& \bl{3} is
\begin{align}
&\sum_{i_{1}+i_{2}+\cdots+i_{r}=N}\frac{a_{1}^{i_{1}}a_{2}^{i_{2}}\cdots a_{r}^{i_{r}}q^{\binom{i_1}{2}+\binom{i_2}{2}+\cdots+\binom{i_r}{2}}}{(q;q)_{i_{1}}(q;q)_{i_{2}}\cdots(q;q)_{i_{r}}}\nonumber\\
&\quad\quad\times \left(1+\sum_{s=1}^{r}q^{i_1+i_2+\cdots+i_s}a_s q^N \prod_{k=1}^{s-1}\left(1+a_k q^N\right)\right).
\end{align}

Finally, we notice that the generating function of $r$-colored strict partitions is
\begin{align}
(-a_{1}q;q)_{\infty}(-a_{2}q;q)_{\infty}\cdots(-a_{r}q;q)_{\infty}.
\end{align}

Hence
\begin{align*}
&(-a_{1}q;q)_{\infty}(-a_{2}q;q)_{\infty}\cdots(-a_{r}q;q)_{\infty}\\
&\quad=1+\sum_{N=1}^{\infty}q^{N^{2}}\prod_{j=1}^{r}(-a_{j}q;q)_{N-1} \sum_{i_{1}+i_{2}+\cdots+i_{r}=N}\frac{a_{1}^{i_{1}}a_{2}^{i_{2}}\cdots a_{r}^{i_{r}}q^{\binom{i_1}{2}+\binom{i_2}{2}+\cdots+\binom{i_r}{2}}}{(q;q)_{i_{1}}(q;q)_{i_{2}}\cdots(q;q)_{i_{r}}}\\
&\quad\quad\quad\quad\quad\times \left(1+\sum_{s=1}^{r}q^{i_1+i_2+\cdots+i_s}a_s q^N \prod_{k=1}^{s-1}\left(1+a_k q^N\right)\right).
\end{align*}

\section{Multi-colored overpartitions}\label{sec:generalization}
In the previous section, the main object we study is $r$-colored strict partitions. We notice that our approach can be naturally adapted to other types of partitions. In particular, if we study multi-colored overpartitions, a more general identity can be deduced.

An \emph{$r$-colored overpartition} means $r$ copies of overpartitions attached with colors $a_{1}$, $a_{2}$, $\ldots$, $a_{r}$. We always assume that only the last occurence of each different part in a different color may be overlined. For instance,
$$\overline{2}_{a_2}+2_{a_1}+\overline{1}_{a_2}+1_{a_1}+\overline{1}_{a_1}$$
is a $2$-colored overpartition of $7$. Here we still assume the following generalized order of parts:
$$1_{a_1}<1_{a_2}<\cdots<1_{a_r}<2_{a_1}<2_{a_2}<\cdots<2_{a_r}<3_{a_1}<\cdots.$$

\begin{figure}[ht]
\caption{Four blocks in an overpartition}\label{fig:block2}
\vspace{1em}
\includegraphics{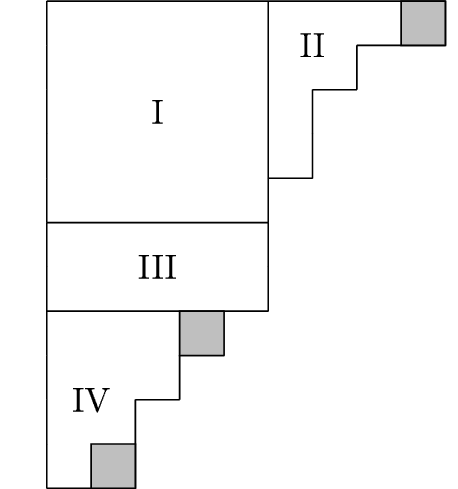}
\end{figure}

We will still use the block decomposition shown in Fig.~\ref{fig:block} as well as the same coloring strategy. To identify the overlined parts, we also shadow the last node of each overlined part in the Ferrers graph (see Fig.~\ref{fig:block2}). Again, we allow $0$ (and hence $\overline{0}$) as a part in Block \bl{2}. In this sense, nodes in Block \bl{1} are neither colored nor shadowed.

Let $N\ge 1$ be the size of the Durfee square $D$, which is also the side length of Block I.

In the following generating functions, for $1\le i\le r$, the exponent of $z_i$ counts the number of overlined parts colored by $a_i$.

\emph{Block \uppercase\expandafter{\romannumeral1}}: From the above arguement, we know that the generating function of $D$ is
\begin{align}
q^{N^{2}}.
\end{align}

\emph{Block \uppercase\expandafter{\romannumeral4}}: It is easy to see that Block \uppercase\expandafter{\romannumeral4} is an $r$-colored overpartition with largest part $\le N-1$. Hence its generating function is
\begin{align}
\prod_{j=1}^{r}\frac{(-a_{j} z_j q;q)_{N-1}}{(a_{j}q;q)_{N-1}}.
\end{align}

\emph{Blocks \uppercase\expandafter{\romannumeral2} \& \uppercase\expandafter{\romannumeral3}}: We start by noticing that the generating function of overpartitions ($0$ not allowed) with at most $i$ parts is (since its conjugate is an overpartition with largest part $\le i$)

\begin{equation}\label{eq:over-1}
\frac{(-zq;q)_i}{(q;q)_i}=\frac{1+z q^i}{1+z}\frac{(-z;q)_i}{(q;q)_i},
\end{equation}
the generating function of overpartitions ($0$ not allowed) with exactly $i$ parts is
\begin{equation}\label{eq:over-2}
\frac{(-zq;q)_i}{(q;q)_i}-\frac{(-zq;q)_{i-1}}{(q;q)_{i-1}}=\frac{q^i (-z;q)_i}{(q;q)_i},
\end{equation}
and the generating function of overpartitions ($0$ allowed) with exactly $i$ parts is
\begin{equation}\label{eq:over-3}
\frac{(1+z)(-zq;q)_{i-1}}{(q;q)_{i-1}}+\frac{q^i (-z;q)_i}{(q;q)_i}=\frac{(-z;q)_i}{(q;q)_i}.
\end{equation}
We have the following two cases:
\begin{enumerate}[\indent 1).]
\item If Block \bl{3} is empty, then thanks to \eqref{eq:over-3}, we know that the generating function of Block \bl{2} is
\begin{align*}
\sum_{i_{1}+i_{2}+\cdots+i_{r}=N}\frac{a_{1}^{i_{1}}a_{2}^{i_{2}}\cdots a_{r}^{i_{r}}(-z_1;q)_{i_{1}}(-z_2;q)_{i_{2}}\cdots(-z_r;q)_{i_{r}}}{(q;q)_{i_{1}}(q;q)_{i_{2}}\cdots(q;q)_{i_{r}}}.
\end{align*}
\item If Block \bl{3} is not empty, then we assume that the part on the top of Block \bl{3} is colored by $a_s$ with $1\le s\le r$. Then the generating function of Block \bl{3} is given by
\begin{align*}
\frac{(1+z_s) a_s q^N}{1-a_s q^N} \prod_{k=1}^{s-1} \frac{1+a_k z_k q^N}{1-a_k q^N}.
\end{align*}
Furthermore, in this case, we only allow $0$ colored by $a_{s}$, $\ldots$, $a_r$ as a part in Block \bl{2} and thoes $0$'s colored by $a_s$, if any, should be non-overlined to ensure that the whole is an $r$-colored overpartition. Suppose there are $i_t$ parts colored by $a_t$ in Block \bl{2} for each $1\le t\le r$, then $i_{1}+i_{2}+\cdots+i_{r}=N$. For $1\le t_1\le s-1$, the overpartition colored by $a_{t_1}$ has exactly $i_{t_1}$ parts and no parts of size $0$, and hence has generating function by \eqref{eq:over-2}
\begin{align*}
\frac{a_{t_1}^{i_{t_1}}q^{i_{t_1}}(-z_{t_1};q)_{i_{t_1}}}{(q;q)_{i_{t_1}}}.
\end{align*}
Next, the overpartition colored by $a_s$ can be treated as an overpartition ($0$ not allowed) with at most $i_s$ parts, and hence has generating function by \eqref{eq:over-1}
\begin{align*}
a_s^{i_s}\frac{1+z_s q^{i_s}}{1+z_s}\frac{(-z_s;q)_{i_s}}{(q;q)_{i_s}}.
\end{align*}
At last, for $s+1\le t_2\le r$, the overpartition colored by $a_{t_2}$ is an overpartition in which we allow $0$ as a part with exactly $i_{t_2}$ parts, and hence has generating function by \eqref{eq:over-3}
\begin{align*}
\frac{a_{t_2}^{i_{t_2}}(-z_{t_2};q)_{i_{t_2}}}{(q;q)_{i_{t_2}}}.
\end{align*}
\end{enumerate}
We conclude that the generating function of Blocks \bl{2} \& \bl{3} is
\begin{align}
&\sum_{i_{1}+i_{2}+\cdots+i_{r}=N}\frac{a_{1}^{i_{1}}a_{2}^{i_{2}}\cdots a_{r}^{i_{r}}(-z_1;q)_{i_{1}}(-z_2;q)_{i_{2}}\cdots(-z_r;q)_{i_{r}}}{(q;q)_{i_{1}}(q;q)_{i_{2}}\cdots(q;q)_{i_{r}}}\nonumber\\
&\quad\quad\times \left(1+\sum_{s=1}^{r}q^{i_1+i_2+\cdots+i_{s-1}}\left(1+z_s q^{i_s}\right) \frac{a_s q^N}{1-a_s q^N}  \prod_{k=1}^{s-1} \frac{1+a_k z_k q^N}{1-a_k q^N}\right).
\end{align}

Since the generating function of $r$-colored overpartitions is
\begin{align}
\frac{(-a_{1} z_1 q;q)_{\infty}(-a_{2} z_2 q;q)_{\infty}\cdots(-a_{r} z_r q;q)_{\infty}}{(a_{1}q;q)_{\infty}(a_{2}q;q)_{\infty}\cdots(a_{r}q;q)_{\infty}},
\end{align}
it follows that \eqref{eq:over-gen} is true and we have completed the proof of Theorem~\ref{thm:CFT}.

\begin{remark}\label{rmk}
The following are special cases of \eqref{eq:over-gen}:

\begin{enumerate}[1).]
\item If we take $z_i=0$ ($1\le i\le r$), then
\begin{align}
&\frac{1}{(a_{1}q;q)_{\infty}(a_{2}q;q)_{\infty}\cdots(a_{r}q;q)_{\infty}}\nonumber\\
&\quad=1+\sum_{N=1}^{\infty}q^{N^{2}}\prod_{j=1}^{r}\frac{1}{(a_{j}q;q)_{N-1}} \sum_{i_{1}+i_{2}+\cdots+i_{r}=N}\frac{a_{1}^{i_{1}}a_{2}^{i_{2}}\cdots a_{r}^{i_{r}}}{(q;q)_{i_{1}}(q;q)_{i_{2}}\cdots(q;q)_{i_{r}}}\nonumber\\
&\quad\quad\quad\quad\quad\times \left(1+\sum_{s=1}^{r}q^{i_1+i_2+\cdots+i_{s-1}}\frac{a_s q^N}{1-a_s q^N} \prod_{k=1}^{s-1}\frac{1}{1-a_k q^N}\right),\label{eq:mul-ptn}
\end{align}
which is a multi-dimensional generalization of Cauchy's identity (cf.~\cite[Eq.~(2.2.8)]{Andr} with $z$ replaced by $aq$):
$$\frac{1}{(aq;q)_\infty}=1+\sum_{N=1}^\infty \frac{a^N q^{N^2}}{(q;q)_N (aq;q)_N}.$$
This multiple summation indeed corresponds to $r$-colored ordinary partitions in our approach.

\item The case $z_i=1$ ($1\le i\le r$) generalizes an identity due to Dousse and Kim (cf.~\cite[Corollary 3.5]{DK}):
\begin{align*}
\dfrac{(-aq;q)_{\infty}}{(aq;q)_{\infty}}=1+\sum_{N=1}^{\infty}\left(\frac{(-q;q)_{N-1}(-aq;q)_{N-1}}{(q;q)_{N-1}(aq;q)_{N-1}}a^{N}q^{N^{2}}+\dfrac{(-q;q)_{N}(-aq;q)_{N}}
{(q;q)_{N}(aq;q)_{N}}a^{N}q^{N^{2}}\right).
\end{align*}
Their proof is based on an overpartition analogue of $q$-binomial coefficients. A further specialization by taking $a=-1$ then recovers \eqref{theta:2}.

\item If we take $a_i=a_i\slash q, z_i=z_iq\:(1\le i\le r)$ and take $q=q^2$, we get the following multi-summation, which can be viewed as the version for {\em ped}, i.e., partitions with even parts distinct.
\begin{align}
&\frac{(-a_{1} z_1 q^2;q^2)_{\infty}(-a_{2} z_2 q^2;q^2)_{\infty}\cdots(-a_{r} z_r q^2;q^2)_{\infty}}{(a_{1}q;q^2)_{\infty}(a_{2}q;q^2)_{\infty}\cdots(a_{r}q;q^2)_{\infty}}\nonumber\\
&\quad=1+\sum_{N=1}^{\infty}q^{2N^{2}-N}\prod_{j=1}^{r}\frac{(-a_{j} z_j q^2;q^2)_{N-1}}{(a_{j}q^2;q^2)_{N-1}} \sum_{i_{1}+i_{2}+\cdots+i_{r}=N}\frac{a_{1}^{i_{1}}\cdots a_{r}^{i_{r}}(-z_1q;q^2)_{i_{1}}\cdots(-z_rq;q^2)_{i_{r}}}{(q^2;q^2)_{i_{1}}\cdots(q^2;q^2)_{i_{r}}}\nonumber\\
&\quad\quad\quad\quad\quad\times \left(1+\sum_{s=1}^{r}q^{2(i_1+i_2+\cdots+i_{s-1})}\left(1+z_s q^{2i_s+1}\right) \frac{a_s q^{2N-1}}{1-a_s q^{2N-1}}  \prod_{k=1}^{s-1} \frac{1+a_k z_k q^{2N}}{1-a_k q^{2N-1}}\right).\label{eq:ped-gen}
\end{align}
Now for the uncolored case $r=1$, we get back to \eqref{theta:4} by setting $a_1=-1,z_1=1$.
\item \eqref{id:alladi} can be deduced from \eqref{eq:over-gen} by taking $a_i \to a_i/z_i$ and then letting $z_i\to \infty$ for $1\le i\le r$.
\end{enumerate}

\end{remark}

\section{Final remarks}

Quite recently, the second and third authors \cite[Eq.~(3.7)]{FT} considered another generalization of Euler's pentagonal number theorem, which involves the numbers of parts and the largest part. On the other hand, we see from \eqref{id:alladi}
\begin{align}
&(-a_{1}yq;q)_{\infty}(-a_{2}yq;q)_{\infty}\cdots(-a_{r}yq;q)_{\infty}\nonumber\\
&\quad=1+\sum_{N=1}^{\infty}\left(y q^N\right)^N\prod_{j=1}^{r}(-a_{j}yq;q)_{N-1} \sum_{i_{1}+i_{2}+\cdots+i_{r}=N}\frac{a_{1}^{i_{1}}a_{2}^{i_{2}}\cdots a_{r}^{i_{r}}q^{\binom{i_1}{2}+\binom{i_2}{2}+\cdots+\binom{i_r}{2}}}{(q;q)_{i_{1}}(q;q)_{i_{2}}\cdots(q;q)_{i_{r}}}\nonumber\\
&\quad\quad\quad\quad\quad\times \left(1+\sum_{s=1}^{r}q^{i_1+i_2+\cdots+i_s}a_s y q^N \prod_{k=1}^{s-1}\left(1+a_k y q^N\right)\right).\label{eq3.1}
\end{align}
Note that \eqref{eq3.1} generalizes \eqref{id:alladi}  in the sense of adding a parameter that counts the number of parts in a partition. We will get back to \eqref{id:alladi} by taking $y=1$. However, it seems to be not easy to consider simultaneously both the number of parts and the largest part, so as to obtain the joint distribution.

At last, it is worth mentioning that in \cite{Nat}, Nataraj established two multivariate generalizations of Euler's pentagonal number theorem related to Rogers--Ramanujan identities.

\subsection*{Acknowledgement}

We would like to acknowledge our gratitude to Ae Ja Yee for her helpful suggestions and comments, which strengthen our original version of Theorem \ref{thm:CFT}. The second and third authors were supported by  National Natural Science Foundation of China grant 11501061.

\end{document}